\documentstyle[12pt, amsfonts]{amsart}
\begin{document}
\def\R{{\mathbb R}}
\def\Z{{\mathbb Z}}
\def\C{{\mathbb C}}
\newcommand{\trace}{\rm trace}
\newcommand{\Ex}{{\mathbb{E}}}
\newcommand{\Prob}{{\mathbb{P}}}
\newcommand{\E}{{\cal E}}
\newcommand{\F}{{\cal F}}
\newtheorem{df}{Definition}
\newtheorem{theorem}{Theorem}
\newtheorem{lemma}{Lemma}
\newtheorem{pr}{Proposition}
\newtheorem{co}{Corollary}
\def\n{\nu}
\def\sign{\mbox{ sign }}
\def\a{\alpha}
\def\N{{\mathbb N}}
\def\A{{\cal A}}
\def\L{{\cal L}}
\def\X{{\cal X}}
\def\F{{\cal F}}
\def\c{\bar{c}}
\def\v{\nu}
\def\d{\delta}
\def\diam{\mbox{\rm dim}}
\def\vol{\mbox{\rm Vol}}
\def\b{\beta}
\def\t{\theta}
\def\l{\lambda}
\def\e{\varepsilon}
\def\colon{{:}\;}
\def\pf{\noindent {\bf Proof :  \  }}
\def\endpf{ \begin{flushright}
$ \Box $ \\
\end{flushright}}

\title[Hyperplane inequality]{A hyperplane inequality for the surface area of projection bodies}

\author{Alexander Koldobsky}

\address{Department of Mathematics\\ 
University of Missouri\\
Columbia, MO 65211}

\email{koldobskiya@@missouri.edu}

\begin{abstract}   We prove a hyperplane inequality for 
the surface area $S(L)$ of a projection body $L$ in $\R^n:$
$$S(L)\ge \frac n{n-1} c_n 
\min_{\xi \in S^{n-1}} S(L\vert \xi^\bot)\  |K|^{1/n},$$
where $|L|$ is the volume of $L,$ $\xi^\bot$ denotes the central 
hyperplane perpendicular to $\xi \in S^{n-1},$ 
$L\vert \xi^\bot$ is the orthogonal projection of $L$ to $\xi^\bot,$ and
$c_n := |B_2^n|^{\frac{n-1}n}/ |B_2^{n-1}| > 1/\sqrt{e}.$ 
The inequality is sharp, with equality when $L=B_2^n$ is the Euclidean ball.
\end{abstract}  
\maketitle

\section{Introduction}

A typical {\it volume comparison problem} 
asks whether inequalities $$f_K(\xi)\le f_L(\xi), \qquad  \forall \xi\in S^{n-1}$$ imply $|K|\le |L|$ 
for any $K,L$ from a certain class of origin-symmetric convex bodies in $\R^n,$ where $f_K$ 
is a certain geometric characteristic of $K.$
One can have in mind the example where $f_K(\xi)= |K\cap \xi^\bot|$ is the hyperplane section function.

In the case where the answer to a volume comparison problem is affirmative,
one can also consider the following {\it separation} problem. Suppose that $\e>0$ and 
\begin{equation}\label{stab2}
f_K(\xi)\le f_L(\xi)-\e, \qquad  \forall \xi\in S^{n-1}.
\end{equation}
Does there exist a constant $c$ not dependent on $\e, K, L$ and such that for every $\e>0$
\begin{equation}\label{volcomp2}
|K|^{\frac{n-1}n} \le |L|^{\frac{n-1}n} - c\e ?
\end{equation}
In the case where the answer is affirmative, assuming that
$$\e= \min_{\xi\in S^{n-1}} \left(f_L(\xi) - f_K(\xi)\right)>0,$$
we get a {\it volume difference inequality}:
\begin{equation}\label{diff2}
 |L|^{\frac{n-1}n} - |K|^{\frac{n-1}n} \ge c\e = c \min_{\xi\in S^{n-1}} \left(f_L(\xi)-f_K(\xi)\right).
\end{equation}
Again, if $f_K$ converges to zero uniformly in $\xi$ when $K$ approaches the empty set,
we get a hyperplane inequality:
\begin{equation}\label{hyper2}
|L|^{\frac{n-1}n} \ge c\min_{\xi\in S^{n-1}} f_L(\xi).
\end{equation}

A separation for the hyperplane projection function
$$f_K(\xi)= P_K(\xi)= \left|K\vert \xi^\perp\right|, \qquad \xi\in S^{n-1},$$
was proved in \cite{K6}.
Here $K\vert\xi^\bot$ is the orthogonal projection of $K$ to $\xi^\bot.$
In Section 2 we reprove this result with the best possible constant. Note that the proof in \cite{K6} essentially provided 
the best possible constant, but along the way the constant was estimated, so the final formulation
was not optimal. Then we apply the corresponding volume difference inequality  to prove 
a hyperplane inequality for the surface
area. We prove that for any projection body $K$ in $\R^n$ the surface area
$$S(K)\ge \frac n{n-1} c_n
\min_{\xi \in S^{n-1}} S(K\vert \xi^\bot)\  |K|^{1/n},$$
with $c_n := |B_2^n|^{\frac{n-1}n}/ |B_2^{n-1}| > 1/\sqrt{e}.$ 
The inequality is sharp, with equality when $K=B_2^n$ is the Euclidean ball.

\section{Main result}

In this section we consider the hyperplane projection function
$$f_K(\xi)= |K\vert\xi^\perp|, $$
where $K\vert\xi^\perp$ is the orthogonal projection of $K$ to the hyperplane $\xi^\perp.$
The corresponding volume comparison result is known as Shephard's problem,
which was posed in 1964 in \cite{Sh} and solved soon after that by Petty \cite{Pe} and 
Schneider \cite{S1}, independently. Suppose that $K$ and $L$ are origin-symmetric convex bodies in 
$\R^n$ so that  $|K\vert \xi^\bot|\le |L\vert \xi^\bot|$ for every $\xi\in S^{n-1}.$ Does it follow that
$|K|\le |L|?$ The answer is affirmative only in dimension 2. Both solutions use the fact that
the answer to Shephard's problem is affirmative in every dimension 
under the additional assumption that $L$ is a projection body; see definition below.

The separation result for projections was proved in \cite{K6}. However, the constant $c_n$
was estimated from below by $1/\sqrt{e}$, so we now formulate the
result with the best possible constant.
\begin{theorem} \label{main-proj1} (\cite{K6}) Suppose that $\e>0$,  $K$ and $L$ are origin-symmetric
convex bodies in $\R^n,$ and $L$ is a projection body.  If for every $\xi\in S^{n-1}$
\begin{eqnarray}\label{proj2}
|K\vert \xi^\bot|\le |L\vert \xi^\bot| - \e,
\end{eqnarray}
then
$$|K|^{\frac{n-1}n}  \le |L|^{\frac{n-1}n} - c_n \e.$$
\end{theorem}

We need several more definitions and results from convex geometry. We refer the reader
to \cite{S2} and \cite{KRZ} for details.

The {\it support function} of a convex body $K$ in $\R^n$ is defined by
$$h_K(x) = \max_{\{\xi\in \R^n:\|\xi\|_K=1\}} (x,\xi),\quad x\in \R^n.$$ 
If $K$ is origin-symmetric, then $h_K$ is a norm on $\R^n.$

The {\it surface area measure} $S(K, \cdot)$ of a convex body $K$ in 
$\R^n$ is defined as follows: for every Borel set $E \subset S^{n-1},$ 
$S(K,E)$ is equal to Lebesgue measure of the part of the boundary of $K$
where normal vectors belong to $E.$ 
We usually consider bodies with absolutely continuous surface area measures.
A convex body $K$ is said to have the {\it curvature function} 
$$ f_K: S^{n-1} \to \R,$$
if its surface area measure $S(K, \cdot)$ is absolutely 
continuous with respect to Lebesgue measure $\sigma_{n-1}$ on 
$S^{n-1}$, and
$$
\frac{d S(K, \cdot)}{d \sigma_{n-1}}=f_K \in L_1(S^{n-1}),
$$
so $f_K$ is the density of $S(K,\cdot).$

By the approximation argument of \cite[Th. 3.3.1]{S2},
we may assume in the formulation of Shephard's problem that the bodies 
$K$ and $L$ are such that  their support functions $h_K,\ h_L$ are 
infinitely smooth functions on $\R^n\setminus \{0\}$.
Using \cite[Lemma 3.16]{K4}
we get in this case that
the Fourier transforms $\widehat{h_K},\ \widehat{h_L}$ are the
extensions of infinitely differentiable functions on the sphere
to homogeneous distributions on $\R^n$ of degree $-n-1.$
Moreover, by a similar approximation argument (see also \cite[Section 5]{GZ}),
we may assume that  our bodies have absolutely continuous surface area 
measures. Therefore, in the rest of this section, $K$ and $L$ are 
convex symmetric bodies with infinitely smooth support functions and absolutely 
continuous surface area measures.

The following version of Parseval's formula was proved in \cite{KRZ} (see also \cite[Lemma 8.8]{K4}):
\begin{equation} \label{pars-proj}
\int_{S^{n-1}} \widehat{h_K} (\xi) \widehat{f_L}(\xi)\ d\xi =
(2\pi)^n \int_{S^{n-1}} h_K(x) f_L(x)\ dx.
\end{equation}

The volume of a body can be expressed in terms of its support function and 
curvature function:
\begin{equation}\label{vol-proj}
|K| = \frac 1n \int_{S^{n-1}}h_K(x) f_K(x)\ dx.
\end{equation}

If $K$ and $L$ are two convex bodies in $\R^n$ the {\it mixed volume} $V_1(K,L)$
is equal to  
$$V_1(K,L)= \frac{1}{n} \lim_{\e\to +0}
\frac{|K+\epsilon L|- |K|}{\e}.$$
We use the 
first Minkowski inequality (see \cite[p.23]{K4}):  
for any convex bodies $K,L$ in $\R^n,$ 
\begin{equation} \label{firstmink}
V_1(K,L) \ge |K|^{(n-1)/n} |L|^{1/n}.
\end{equation}
The mixed volume can be expressed in terms of the support and
curvature functions:

\begin{equation}\label{mixvol-proj}
V_1(K,L) = \frac 1n \int_{S^{n-1}}h_L(x) f_K(x)\ dx.
\end{equation}

Let $K$ be an origin-symmetric convex body in $\R^n.$ The {\it
projection body} $\Pi K$ of $K$ is defined as an origin-symmetric convex 
body in $\R^n$ whose support function in every direction is equal to
the volume of the hyperplane projection of $K$ to this direction, i.e.
for every $\theta\in S^{n-1},$ 
\begin{equation} \label{def:proj}
h_{\Pi K}(\theta) = |K\vert\theta^{\perp}|.
\end{equation}
If $L$ is the projection body of some convex body, we simply say 
that $L$ is a projection body.  The Minkowski (vector) sum of projection bodies
is also a projection body; see for example \cite[p. 149]{G3}. 
\medbreak
\noindent {\bf Proof of Theorem \ref{main-proj1}.} It was proved in \cite{KRZ} that
\begin{equation} \label{f-proj}
P_K(\xi)=|K\vert \xi^\bot| = -\frac 1{\pi} \widehat{f_K}(\xi),\qquad \forall \xi\in S^{n-1},
\end{equation}
where $f_K$ is extended from the sphere to a homogeneous function of degree 
$-n-1$ on the whole $\R^n,$ and the Fourier transform $\widehat{f_K}$ is the 
extension of a continuous function $P_K$ on the sphere to a homogeneous of degree 1
function on $\R^n.$

Therefore, the condition (\ref{proj2}) can be written as
\begin{equation} \label{fourier-proj}
\frac 1{\pi}\ \widehat{f_K}(\xi) \ge  \frac 1{\pi}\ \widehat{f_L}(\xi) + \e, \qquad \forall \xi\in S^{n-1}.
\end{equation}
It was also proved in \cite{KRZ} that an infinitely smooth origin-symmetric convex body 
$L$ in $\R^n$ is a projection body if and only if 
$\widehat{h_L} \le 0$ on the sphere $S^{n-1}.$ Integrating (\ref{fourier-proj})
with respect to a negative density we get
$$\int_{S^{n-1}} \widehat{h_L}(\xi) \widehat{f_L}(\xi)\ d\xi \ge \int_{S^{n-1}} \widehat{h_L}(\xi) \widehat{f_K}(\xi)\ d\xi 
- \pi\e \int_{S^{n-1}} \widehat{h_L}(\xi)\ d\xi.$$
Using this, (\ref{vol-proj}) and (\ref{pars-proj}), we get
$$ (2\pi)^n n |L| = (2\pi)^n \int_{S^{n-1}} h_L(x) f_L(x)\  dx =
\int_{S^{n-1}} \widehat{h_L}(\xi) \widehat{f_L}(\xi)\ d\xi$$
$$\ge \int_{S^{n-1}} \widehat{h_L}(\xi) \widehat{f_K}(\xi)\ d\xi  - \pi\e\int_{S^{n-1}}\widehat{h_L}(\xi)\ d\xi$$
\begin{equation} \label{eq31}
=(2\pi)^n \int_{S^{n-1}} h_L(x) f_K(x)\  dx - \pi\e\int_{S^{n-1}}\widehat{h_L}(\xi)\ d\xi.
\end{equation}
We estimate the first summand from below using the first Minkowski inequality:
\begin{equation} \label{eq32}
(2\pi)^n \int_{S^{n-1}} h_L(x) f_K(x)\  dx \ge (2\pi)^n n \left(\vol_n(L)\right)^{\frac 1n} \left(\vol_n(K)\right)^{\frac {n-1}n}.
\end{equation}

To estimate the second term in (\ref{eq31}), note that, by (\ref{f-proj}),
the Fourier transform of the curvature function of the Euclidean ball
$$\widehat{f_2}(\xi) = -\pi |B_2^{n-1}|,\qquad \forall \xi\in S^{n-1}.$$
Therefore, by Parseval's formula, (\ref{mixvol-proj})  and the first Minkowski inequality,
$$\pi \e \int_{S^{n-1}}\widehat{h_L}(\xi)\ d\xi = - \frac {\e}{|B_2^{n-1}|}
\int_{S^{n-1}}\widehat{h_L}(\xi)\widehat{f_2}(\xi)\ d\xi$$
$$= -  \frac {(2\pi)^n \e}{|B_2^{n-1}|}
\int_{S^{n-1}} h_L(x) f_2(x)\ dx
= -  \frac {(2\pi)^n \e}{|B_2^{n-1}|} nV_1(B_2^n,L)$$
$$\le  - \frac {(2\pi)^n n\e}{|B_2^{n-1}|} |L|^{\frac 1n}
 |B_2^n|^{\frac {n-1}n}
 = - (2\pi)^n n \e c_n |L|^{\frac 1n}.$$
 Combining this with (\ref{eq31}) and (\ref{eq32}), we get the result.
 \endpf
 
As explained in the Introduction, the separation result of Theorem \ref{main-proj1} leads to
a volume difference inequality of the type (\ref{diff2}). 

\begin{co} If $L$ is a projection body in $\R^n$ and $K$ is an arbitrary origin-symmeric convex body in $\R^n$
so that $$\min_{\xi\in S^{n-1}} (|L\vert \xi^\bot| - |K\vert \xi^\bot|) >0,$$ then
\begin{equation} \label{diff-proj}
|L|^{\frac {n-1}n} - |K|^{\frac {n-1}n} \ge c_n \min_{\xi\in S^{n-1}} (|L\vert \xi^\bot| - |K\vert \xi^\bot|).
\end{equation}
\end{co}
Sending $K$ to the empty set in (\ref{diff-proj}), we get a hyperplane inequality of the type (\ref{hyper2})
which was earlier deduced directly from the affirmative part of Shephard's problem in \cite[Corollary 9.3.4]{G3}:
if $L$ is a projection body in $\R^n$, then
\begin{equation} \label{hyper-proj}
|L|^{\frac {n-1}n} \ge c_n \min_{\xi\in S^{n-1}} |L\vert \xi^\bot|.
\end{equation}
Recall that $c_n > 1/\sqrt{e}.$ For general symmetric convex bodies, $c_n$ can be replaced in (\ref{hyper-proj}) by
 $c_n/n^{(n-1)/2n}$ (see \cite[Remark 9.3.5]{G3}; note that
 with an extra factor $(3/2)^{(n-1)/n}$ in the left-hand side the estimate follows from a result of Ball \cite[p. 899]{Ba2}).  
 Moreover, Ball \cite[Theorem 5]{Ba2} proved that the constant of the order $1/\sqrt{n}$ is optimal for general 
symmetric convex bodies, namely
there exists a constant $\delta>0$ such that for every $n\in \N$ there is an origin-symmetric 
convex body $L$ in $\R^n$ such that
$$ |L|^{\frac {n-1}n} \le \frac 1{\delta \sqrt{n}} \min_{\xi\in S^{n-1}} |L\vert \xi^\bot|.$$

The volume difference inequality (\ref{diff-proj}) allows to prove a hyperplane
inequality for the surface area of projection bodies.
\begin{co} \label{surf-proj}Let $L$ be a projection body in $\R^n,$ then the surface area
$$S(L) \ge \frac n{n-1} c_n \min_{\xi\in S^{n-1}} S(L\vert\xi^\bot)|L|^{\frac 1n}.$$
\end{co}
\pf The surface area of $L$ can be computed as
$$S(L) = \lim_{\e \to +0} \frac {\left|L+\e B_2^n\right| - \left|L\right|}{\e}.$$
For every $\e>0$ the Minkowski sum $L+\e B_2^n$ is also a projection body. 
The inequality (\ref{diff-proj}) with the bodies
$L+\e B_2^n$ and $L$ in place of $L$ and $K$ implies
$$\frac {|L+\e B_2^n|^{\frac {n-1}n} - |L|^{\frac {n-1}n}}{\e} \ge
c_n \min_{\xi\in S^{n-1}} \frac {|(L\vert \xi^\bot)+\e B_2^{n-1}| - 
|L\vert \xi^\bot|}{\e}.$$
By the Minkowski theorem on mixed volumes (\cite[Theorem 5.1.6]{S2} or \cite[Theorem A.3.1]{G3}),
\begin{equation} \label{quer}
\frac{|(L\vert \xi^\bot)+ \e B_2^{n-1}|-|L\vert \xi^\bot|}{\e} = 
\sum_{i=1}^{n-1} {n-1 \choose i} W_i(L\vert \xi^\bot) \e^{i-1},
\end{equation}
where $W_i$ are quermassintegrals. The function $\xi\mapsto L\vert \xi^\bot$ is continuous 
from $S^{n-1}$ to the class of origin-symmetric convex sets equiped with the Hausdorff metric,
and $W_i$'s are also continuous with respect to this metric (see \cite[p. 275]{S2}), so the 
functions $\xi\mapsto W_i(L\vert \xi^\bot)$ are continuous and, hence, bounded on the sphere.
This implies that the left-hand side of (\ref{quer}) converges to $S(L\vert \xi^\bot)$ as $\e\to 0$
uniformly with respect to $\xi.$ The latter allows to switch the limit and maximum  
in the right-hand side of (\ref{surflimit}), as $\e\to 0$.
Sending $\e$ to zero in (\ref{surflimit}), we get
$$\frac {n-1}n |L|^{-1/n} S(L) \ge c_n  \min_{\xi \in S^{n-1}} S(L\vert \xi^\bot).$$
\endpf

\bigbreak

{\bf Acknowledgement.} This work was supported in part by
the US National Science Foundation, grant DMS-1001234.

\end{document}